\documentclass{article}\usepackage{amsmath,amssymb,latexsym,amsfonts,enumerate,theorem,longtable,amscd}\newtheorem{definition}{Definition} 
\newtheorem{theorem}{Theorem}\newtheorem{proposition}{Proposition}\newtheorem{corollary} {Corollary}
\newtheorem{lemma} {Lemma}
\begin{document}

\begin{center}
{\bf\large Covariant Differential Calculi on $SP_q^{1|2}$}
\end{center}

\begin{center}
Salih Celik

Department of Mathematics, Yildiz Technical University, DAVUTPASA-Esenler, Istanbul, 34210 TURKEY.
\end{center}

\noindent{\bf MSC:} 17B37; 81R60

\noindent{\bf Keywords:} Quantum symplectic superspace, super $\star$-algebra, differential calculus, quantum supergroup $SP_q(1|2)$, unitary orthosymplectic quantum supergroup

\begin{abstract}
A unitary orthosymplectic quantum supergroup is introduced. Two covariant differential calculi on the quantum superspace $SP_q^{1|2}$ are presented. The $h$-deformed symplectic superspaces via a contraction of the $q$-deformed symplectic superspaces are obtained. A new $h$-deformation of the Heisenberg superalgebra is given.
\end{abstract}

\section{Introduction}\label{intro}

The theory of compact matrix quantum groups developed by Woronowicz \cite{Woronowicz1}. Quantum (super)groups have a rich mathematical structure. Since they may be regarded as a noncommutative extension of Lie (super)groups they provide a very powerful tool for investigations of noncommutative geometry. Due to pioneering work by Woronowicz  \cite{Woronowicz2}, the quantum (super)groups supplied concrete examples of noncommutative differential geometry \cite{Connes} by introducing a consistent differential calculus on the noncommutative spaces of the quantum (super)groups. In this approach the quantum group is taken as the basic noncommutative space and the differential calculus on the group is deduced from its properties.

A quantum space is a space that quantum group acts with linear transformations and whose coordinates belong to a noncommutative associative algebra \cite{Manin1}. Interpreting the dual space of the quantum space as differentials of the coordinates, covariant differential calculus on the quantum space has been developed in \cite{Wess-Zumino}. The natural extension of their scheme to superspace \cite{Manin2} was introduced in \cite{Soni}.

The quantum superplane is the simplest example of a noncommutative superspace. The noncommutative geometry of the quantum superplane investigated in \cite{Celik1}. In this paper, we have investigated the noncommutative geometry of the quantum symplectic $(1+2)$-superspace, denoted by $SP_q^{1|2}$ and given some related topics.

\section{Review of Quantum Symplectic Group}\label{sec:2}

In this section, we will give some information about the structures of classical and quantum symplectic groups as much as needed.

\subsection{Classical Case}\label{sec2.1}

In the classic case, the algebra ${\cal O}(SP(1|2))$ consists of five even and four odd generators. If even generators are denoted by $a,b,c,d,e$ and odd generators by $\alpha,\beta,\gamma,\delta$, an element of the symplectic group $SP(1|2)$ can be written as
\begin{equation} \label{1}
T = \begin{pmatrix} a & \alpha & b \\ \gamma & e & \beta \\ c & \delta & d \end{pmatrix} =(t_{ij}).
\end{equation}
The orthosymplectic condition is
\begin{equation} \label{2}
T^{st}CT = DC, \quad TC^{-1}T^{st} = DC^{-1}
\end{equation}
where $T^{st}$ denotes the super transposition of $T$. Here, the superdeterminant of $T$ is defined by
\begin{equation*}
D = ad - bc - \alpha\delta = da - cb + \delta\alpha = {\bf 1}.
\end{equation*}
Therefore, in a matrix belonging to orthosymplectic group $OSP(1|2)$, there are six elements, four of which are even and other two are odd:
\begin{equation*}
\gamma = a\delta - c \alpha, \quad e = {\bf 1}-\alpha\delta, \quad \beta = b\delta - d\alpha.
\end{equation*}
The matrix $C$ of metric is given by
\begin{equation} \label{3}
C = \begin{pmatrix} 0 & 0 & -1 \\ 0 & 1 & 0 \\ 1 & 0 & 0 \end{pmatrix}.
\end{equation}
We refer to the basic papers \cite{Ber-Tos}, \cite{Fra-Sci-Sor}, \cite{Rit-Sch} about the supergroup $OSP(1|2)$ for interested readers.

\subsection{Quantum Case}\label{sec2.2}

The algebra ${\cal O}(OSP_q(1|2))$ is generated by the even elements $a,b,c,d$ and odd elements $\alpha,\delta$. Standard FRT construction \cite{FRT} is
obtained via the matrix $R$ given in \cite{Kul-Res}. Using the RTT-relations and the $q$-orthosymplectic contition, all defining relations of
${\cal O}(OSP_q(1|2))$ are explicitly obtained in \cite{Aiz-Chak1}:
\begin{theorem} 
The generators of ${\cal O}(SP_q(1|2))$  satisfy the relations
\begin{eqnarray} \label{4}
a b &=& q^2 b a, \quad a c = q^2 c a, \quad a \alpha = q \alpha a, \quad a \delta = q \delta a + (q-q^{-1}) \alpha c, \nonumber\\
ad &=& da + (q-q^{-1})[(1+q^{-1}) bc + q^{-1/2} \alpha\delta], \nonumber\\
b c &=& c b, \quad bd = q^2 db, \quad b \alpha = q^{-1} \alpha b, \quad b \delta = q \delta b, \nonumber\\
cd &=& q^2 dc, \quad c \alpha = q^{-1} \alpha c, \quad c \delta = q \delta c, \nonumber\\
d \alpha &=& q^{-1} \alpha d + (q^{-1}-q) \delta b, \quad d \delta = q^{-1} \delta d, \\
\alpha\delta &=& -q \delta\alpha + q^{-1/2}(q-q^{-1}) bc, \nonumber\\
\alpha^2 &=& q^{-1/2}(q-1) ba, \quad \delta^2 = q^{-1/2}(q-1) dc.\nonumber
\end{eqnarray}
\end{theorem}

In (\ref{4}), the relations involving the elements $\gamma$, $e$ and $\beta$ are not written. They and identities from quantum analogue of (\ref{2}) with (\ref{3}) can be found in \cite{Aiz-Chak1}.

\noindent{\bf Note 1.} The relations (\ref{4}) are connected to the isomorphism interchanging $a$ and $d$, $b$ and $c$, $\alpha$ and $\delta$, $q$ and $q^{-1}$.

The quantum superdeterminant is defined by
\begin{equation*}
D_q = ad - qbc - q^{1/2}\alpha\delta = da - q^{-1}bc + q^{-1/2}\delta\alpha.
\end{equation*}
This element of ${\cal O}(SP_q(1|2))$ commutes with all elements of ${\cal O}(SP_q(1|2))$.

If ${\cal A}$ and ${\cal B}$ are $Z_2$-graded algebras, then their tensor product ${\cal A}\otimes{\cal B}$ is the $Z_2$-graded algebra whose underlying space is $Z_2$-graded tensor product of ${\cal A}$ and ${\cal B}$. The following definition gives the product rule for tensor product of algebras. Let us denote by $p(a)$ the {\it grade} (or {\it degree}) of an element $a\in {\cal A}$.

\begin{definition}
If ${\cal A}$ is a $Z_2$-graded algebra, then the product rule in the $Z_2$-graded algebra ${\cal A}\otimes {\cal A}$ is defined by
\begin{equation*}
(a_1\otimes a_2)(a_3\otimes a_4) = (-1)^{p(a_2)p(a_3)} a_1a_3\otimes a_2a_4
\end{equation*}
where $a_i$'s are homogeneous elements in the algebra ${\cal A}$.
\end{definition}

The super-Hopf algebra structure of ${\cal O}(SP_q(1|2))$ is given as usual in the following proposition.
\begin{theorem} 
There exists a unique super-Hopf algebra structure on the superalgebra ${\cal O}(SP_q(1|2))$ with co-maps $\Delta$, $\epsilon$ and $\kappa$ such that
\begin{eqnarray*}
\Delta(t_{ij}) = \sum_{k=1}^3 t_{ik} \otimes t_{kj}, \quad \epsilon(t_{ij}) = \delta_{ij}, \quad S(T)=T^{-1}
\end{eqnarray*}
where
\begin{equation} \label{5}
T^{-1} = \begin{pmatrix} d & q^{-1/2}\beta & -q^{-1}b \\ -q^{1/2}\delta & e & q^{-1/2}\alpha \\ -qc & -q^{1/2}\gamma & a \end{pmatrix}.
\end{equation}
\end{theorem}

\section{Quantum Symplectic Superspaces}\label{sec3}

Some spaces do not agree the quantum deformation for each value of the deformation parameter $q$. In other words, some deformations such as quantum symplectic superspace may exist for special choice of $q$. In this section, one has defined the function algebras of the quantum symplectic superspace and its dual.

The elements of the symplectic superspace are supervectors generated by an even and two odd components. One defined a Z$_2$-graded symplectic space $SP^{1|2}$ by dividing the superspace $SP^{1|2}$ of 3x1 matrices into two parts $SP^{1|2}=V_0\oplus V_1$. A vector is an element of $V_0$ (resp. $V_1$) and is of grade 0 (resp. 1) if it has the form shown below:
\begin{equation*}
\begin{pmatrix} 0 \\ x \\ 0\end{pmatrix}, \quad \mbox{resp.} \quad \begin{pmatrix} \xi \\ 0 \\ \eta \end{pmatrix}.
\end{equation*}
While the even element $x$ commutes to other two elements, the odd elements satisfy the relations $\Theta_i\Theta_j=-\Theta_j\Theta_i$ for $i,j=1,2$.
Let us denote the exterior algebra of the space $SP^{1|2}$ by $\Lambda(SP^{1|2})$. Then, $\Lambda(SP^{1|2})$ is a commutative superalgebra generated by the odd generator $\theta$ and even generators $y, z$ where $\theta^2=0$.

\subsection{The Algebra of Polynomials on the Quantum Superspace $SP_q^{1|2}$}\label{sec3.1}

\begin{definition}[\cite{Aiz-Chak1}] 
Let $k\{\xi,x,\eta\}$ be a free associative algebra generated by $x$, $\xi$, $\eta$ and $I_q$ is a two-sided ideal generated by $x\xi-q\xi x$,
$x\eta-q^{-1}\eta x$, $\xi\eta+q^{-2}\eta\xi+q^{-2}Qx^2$, $\xi^2$ and $\eta^2$. The quantum superspace $SP_q^{1|2}$ with the function algebra
$${\cal O}(SP_q^{1|2}) = k\{\xi,x,\eta\}/I_q$$
is called Z$_2$-graded quantum symplectic space (or quantum symplectic superspace) where  $Q=q^{1/2}-q^{3/2}$ and $q\ne0$.
\end{definition}

According to Definition 2, if $(\xi,x,\eta)^t\in SP_q^{1|2}$ then we have
\begin{equation} \label{6}
x\xi = q \xi x, \quad x\eta = q^{-1} \eta x, \quad \eta\xi = - q^2 \xi\eta + q^{1/2}(1-q) x^2, \quad \xi^2 = 0 = \eta^2.
\end{equation}
This associative algebra over the complex number is known as the algebra of polynomials over quantum symplectic (1+2)-superspace.

It is easy to see the existence of representations that satisfy (\ref{6}); for instance, the three 2x2 matrix representations of the coordinate functions $\xi$, $x$ and $\eta$
\begin{equation*}
\rho(\xi) = \begin{pmatrix} 0 & q^{1/2} \\ 0 & 0 \end{pmatrix}, \quad \rho(x)=\begin{pmatrix} q & 0 \\ 0 & 1 \end{pmatrix}, \quad
\rho(\eta) = \begin{pmatrix} 0 & 0 \\ 1-q & 0 \end{pmatrix}
\end{equation*}
satisfy the relations (\ref{6}). Similarly, the three 3x3 matrix representations
\begin{equation*}
\rho(\xi) = \begin{pmatrix} 0 & 0 & q^{1/2} \\ 0 & 0 & 0 \\ 0 & 0 & 0 \end{pmatrix}, \quad \rho(x)=\begin{pmatrix} q & 0 & 0 \\ 0 & 0 & 0 \\ 0 & 0 & 1 \end{pmatrix}, \quad
\rho(\eta) = \begin{pmatrix} 0 & 0 & 0\\ 0 & 0 & 0\\ 1-q & 0 & 0 \end{pmatrix}
\end{equation*}
also satisfy the relations (\ref{6}).

The set $\{x^k\xi^l\eta^m: \, k\in{\mathbb N}_0, l,m=0,1\}$ forms a vector space basis of ${\cal O}(SP_q^{1|2})$.

\noindent{\bf Note 2.} The relations (\ref{6}) are invariant if we replace $q$ by $q^{-1}$, $x$ by $\sqrt{-1} \,x$ and $\xi$ by $\eta$ in the relations (\ref{6}) which are compatible with  Note 1.

Consideration of quantum group extensions of some models of integrable quantum field theories with $OSP(1|2)$ symmetries where superspheres
appear \cite{Saleur-Kaufmann} will require quantum superspheres. The quantum symplectic supersphere can be defined by using the quantum version of the matrix $C$ in (\ref{3}),
\begin{equation} \label{7}
C_q = \begin{pmatrix} 0 & 0 & -q^{-1/2} \\ 0 & 1 & 0 \\ q^{1/2} & 0 & 0 \end{pmatrix}.
\end{equation}

\begin{definition}[\cite{Aiz-Chak2}] 
The quantum supersphere on the quantum symplectic superspace is defined by
\begin{equation*}
r = X^{st}C_qX = q^{-1/2} \xi\eta + x^2 - q^{1/2} \eta\xi.
\end{equation*}
\end{definition}

The definition of dual quantum symplectic superspace is as follows.
\begin{definition}
$\Lambda(SP_q^{2|1})$ is the superalgebra with odd generator $\theta$ and even generators $y$, $z$ and the following set of relations
\begin{equation} \label{8}
\theta y = q^{-1} y\theta, \quad \theta z = q z\theta, \quad yz = q^2 zy, \quad \theta^2 = q^{1/2}(q-1) zy.\end{equation}
We call $\Lambda(SP_q^{2|1})$ the exterior superalgebra of the quantum symplectic superspace $SP_q^{1|2}$
\end{definition}

\noindent{\bf Note 3.} The relations (\ref{8}) are invariant if we replace $q$ by $q^{-1}$, $\theta$ by $\sqrt{-1} \,\theta$ and $y$ by $z$ in the relations (\ref{11}).

\subsection{Coactions on the Quantum Symplectic Superspace}\label{sec3.2}

Let $a,b,c,d,e,\gamma,\alpha,\delta,\beta$ be elements of an algebra ${\cal A}$. Assuming that the generators of ${\cal O}(OSP_q(1|2))$ super-commute with the elements of ${\cal O}(SP_q^{1|2})$ define $\xi',x',\eta'$ and $\xi'',x'',\eta''$ using the following matrix equalities
\begin{equation} \label{9}
\begin{pmatrix} \xi' \\ x' \\ \eta' \end{pmatrix}=\begin{pmatrix} a & \alpha & b \\ \gamma & e & \beta \\ c & \delta & d \end{pmatrix} \begin{pmatrix} \xi\\ x \\ \eta \end{pmatrix} \,\, \mbox{and} \,\, (\xi'', x'', \eta'') = (\xi, x, \eta) \begin{pmatrix} a & \alpha & b \\ \gamma & e & \beta \\ c & \delta & d \end{pmatrix}.
\end{equation}
If we assume that $q\ne1$ then we have the following proposition proving straightforward computations as follows.

\begin{theorem} 
If the couples $(\xi', x', \eta')$ and $(\xi'', x'', \eta'')$ in (\ref{9}) satisfy the relations (\ref{6}), then the generators of ${\cal O}(SP_q(1|2))$ fulfill the relations (\ref{4}).
\end{theorem}

A left quantum space for a Hopf algebra $H$ is an algebra ${\cal X}$ together with an algebra homomorphism (left coaction)
$\delta_L:{\cal X}\longrightarrow H\otimes{\cal X}$ such that
\begin{equation*}
(\mbox{id}\otimes\delta_L)\circ\delta_L=(\Delta\otimes\mbox{id})\circ\delta_L \quad \mbox{and} \quad (\epsilon\otimes\mbox{id})\circ\delta_L=\mbox{id}.
\end{equation*}

\begin{theorem} 
(i) The algebra ${\cal O}(SP_q^{1|2})$ is a left and right comodule algebra of the Hopf algebra ${\cal O}(SP_q(1|2))$ with left coaction
$\delta_L$ and right coaction $\delta_R$ such that
\begin{equation} \label{10}
\delta_L(X_i) = \sum_{k=1}^3 t_{ik} \otimes X_k, \quad \delta_R(X_i) = \sum_{k=1}^3 X_k \otimes t_{ki}.
\end{equation}

\noindent(ii) The quantum supersphere $r$ belongs to the center of the algebra ${\cal O}(SP_q^{1|2})$ and satisfies $\delta_L(r)={\bf 1}\otimes r$ and $\delta_R(r)=r \otimes{\bf 1}$.
\end{theorem}

\noindent{\it Proof}
(i): These assertions are obtained from the relations in (\ref{4}) with (\ref{6}).

(ii): That $r$ is a central element of ${\cal O}(SP_q^{1|2})$ which is shown by using the relations in (\ref{6}). To show that $\delta_L(r)={\bf 1}\otimes r$ and $\delta_R(r)=r \otimes{\bf 1}$ we use the definitions of $\delta_L$ and $\delta_R$ in (\ref{10}) and the relations (\ref{4}) with $D_q={\bf 1}$.

\section{Unitary Orthosymplectic Quantum \\ Supergroup}\label{sec4}

In this section, by defining a $Z_2$-graded involution on the superalgebra \\${\cal O}(OSP_q(1|2))$, we will obtain unitary orthosymplectic quantum supergroup.

It is possible to define the involution on the Grassmann generators. However, there are two possibilities to do so. The first is given below. Other will be
given in section 6.

There exist two cases according to the values of the deformation parameter $q$. We assume that here the deformation parameter $q$ is a real number. Let
${\cal A}$ be a $Z_2$-graded algebra.

\begin{definition} 
A conjugate-linear map $\vartheta \mapsto \vartheta^\star$ of degree zero is called a $Z_2$-graded involution on $Z_2$-graded algebra ${\cal A}$ if
\begin{equation*}
(\vartheta_1\vartheta_2)^\star = (-1)^{p(\vartheta_1)p(\vartheta_2)} \vartheta_2^\star \vartheta_1^\star, \quad (\vartheta^\star)^\star = (-1)^{p(\vartheta)}\vartheta
\end{equation*}
for all $\vartheta,\vartheta_1, \vartheta_2 \in {\cal A}$. The pair $({\cal A},\star)$ is called a $Z_2$-graded $\star$-algebra.
\end{definition}
Let ${\cal A}$ is a $Z_2$-graded Hopf algebra and $\star:{\cal A} \longrightarrow {\cal A}$ is an involution such that $({\cal A},\star)$ becomes a $Z_2$-graded Hopf $\star$-algebra. (see, Definition 6) A $Z_2$-graded quantum space ${\cal X}$ for a $Z_2$-graded Hopf $\star$-algebra $\cal A$ is called a $Z_2$-graded $\star$-quantum space if $\cal X$ is a $\star$-algebra and the coaction $\delta_.$ of ${\cal A}$ on ${\cal X}$ satisfies $\delta_.(x^\star)=(\delta_.(x))^\star$ for all $x\in {\cal X}$.

\begin{proposition} 
If $q>0$ then the algebra ${\cal O}(SP_q^{1|2})$ supplied with the $Z_2$-graded involutions determined by
\begin{equation} \label{11}
\xi^\star = q^{1/2}\eta, \quad x^\star = x, \quad \eta^\star = -q^{-1/2}\xi
\end{equation}
becomes a super $\star$-algebra.
\end{proposition}

\noindent{\it Proof}
We must show that the relations (\ref{6}) are invariant under these involutions. By Definition 4.1, since $(\xi\eta)^\star=-\eta^\star\xi^\star=\xi\eta$ and $(\eta\xi)^\star=-\xi^\star\eta^\star=\eta\xi$, third relation in (\ref{6}) is indeed invariant under the star operation. Other relations in (\ref{6}) are also provided in a similar manner.

The supergroup $UOSP(1|2)$ in contrast is made of complex super transformations satisfying
\begin{equation*}
T^{st}CT = C, \quad TT^\ddag = I
\end{equation*}
where $T^\ddag = (T^{st})^\star$.

Let $q>0$. Then the super-Hopf algebra ${\cal O}(SP_q(1|2))$ can be supplied with involutions and it becomes a super-Hopf $\star$-algebra. The corresponding super-Hopf $\star$-algebra is denoted by ${\cal O}(UOSP_q(1|2))$. The matrix identity in the following lemma is obtained from the condition
$TT^\ddag = I$ with (\ref{5}).

\begin{lemma} 
If $T \in OSP_q(1|2)$ then we have
\begin{equation*}
T^\ddag = \begin{pmatrix} a^\star & \gamma^\star & c^\star \\ -\alpha^\star & e^\star & -\delta^\star \\ b^\star & \beta^\star & d^\star \end{pmatrix} = \begin{pmatrix} d & q^{-1/2}\beta & -q^{-1}b \\ -q^{1/2}\delta & e & q^{-1/2}\alpha \\ -qc & -q^{1/2}\gamma & a \end{pmatrix}
\end{equation*}
where $q>0$.
\end{lemma}

According to the lemma above, for a matrix $T$ in $UOSP_q(1|2)$ we can write
\begin{equation*}
T = \begin{pmatrix} a & \alpha & b \\ -q^{-1/2}\beta^\star & e & \beta \\ -q^{-1}b^\star & q^{-1/2}\alpha^\star & a^\star \end{pmatrix}.
\end{equation*}
In this case, the quantum superdeterminant becomes
\begin{equation*}
D_q(T) = aa^\star + bb^\star - \alpha \alpha^\star = {\bf 1} = a^\star a + q^{-2}b^\star b + q^{-1}\alpha^\star \alpha.
\end{equation*}
Eliminating $\alpha^\star \alpha$ from $D_q(T)$ results in
\begin{equation*}
aa^\star - q^2 a^\star a = (1-q^2) ({\bf 1} + b^\star b).
\end{equation*}
For real $q$ such that $0<q<1$  rescaling of $a$ and $b$ by
\begin{equation*}
a \mapsto (1-q^2)^{1/2} \,a, \quad b \mapsto (1-q^2)^{1/2} \,b
\end{equation*}
will give
\begin{equation*}
aa^\star - q^2 a^\star a = {\bf 1} + (1-q^2) b^\star b.
\end{equation*}
This is the form of the two dimensional $q$-oscillator which has been conventionally used in various works.

\begin{definition} 
If ${\cal A}$ is a $\star$-algebra then the involution of ${\cal A}\otimes{\cal A}$ is defined by $(a\otimes b)^\star=a^\star\otimes b^\star$ for all elements $a,b \in {\cal A}$.
A Hopf algebra ${\cal A}$ is called a $\star$-Hopf algebra if ${\cal A}$ is supplied with an involution $\star$ such that
$\Delta(a^\star)=[\Delta(a)]^\star$ and $\epsilon(a^\star)=\overline{\epsilon(a)}$ for $a \in {\cal A}$.
\end{definition}
Specially, in any super-Hopf $\star$-algebra ${\cal A}$, we have $S(S(a^\star)^\star) = (-1)^{p(a)} a$ for $a \in {\cal A}$.

\begin{proposition} 
(i) The generators of ${\cal O}(UOSP_q(1|2))$  satisfy the relations
\begin{eqnarray} \label{12}
a b &=& q^2 b a, \quad a \alpha = q \alpha a, \quad b \alpha = q^{-1} \alpha b, \nonumber\\
a b^\star &=& q^2 b^\star a, \quad b b^\star = b^\star b, \quad b \alpha^\star = q \alpha^\star b, \nonumber\\
a \alpha^\star &=& q \alpha^\star a + q^{1/2}(q^{-1}-q) b^\star \alpha, \quad \alpha \alpha^\star = -q \alpha^\star \alpha + (q^{-1}-q) b^\star b, \nonumber\\
aa^\star &=& a^\star a + (q^{-1}-q) \left[(1+q^{-1})b^\star b + \alpha^\star \alpha\right] \\
\alpha^2 &=& q^{-1/2}(q-1) ba, \quad (\alpha^\star)^2 = q^{-1/2}(1-q) a^\star b^\star.\nonumber
\end{eqnarray}

(ii) The algebra ${\cal O}(UOSP_q(1|2))$ is a super-Hopf $\star$-algebra.
\end{proposition}

\noindent{\it Proof}
(i): The relations (\ref{12}) can be obtained using the relations (\ref{4}).

(ii): To show that ${\cal O}(UOSP_q(1|2))$ is a $\star$-algebra, one must check the demands of Definition 3.6. For examples,
\begin{eqnarray*}
(\alpha^\star)^\star = (q^{1/2}\delta)^\star = q^{1/2} (-q^{-1/2} \alpha) = -\alpha, \\
(b^\star)^\star = (-qc)^\star = -q (-q^{-1} b) = b.
\end{eqnarray*}
All relations in (\ref{12}) are invariant under the $\star$ operation. ${\cal O}(UOSP_q(1|2))$ is obviously a super-Hopf algebra. On the other hand, for example, since
\begin{equation*}
\Delta(a) = a\otimes a - q^{-1/2}\alpha\otimes\beta^\star - q^{-1}b\otimes b^\star,
\end{equation*}
we have
\begin{equation*}
[\Delta(a)]^\star = a^\star\otimes a^\star + q^{-1/2}\alpha^\star\otimes\beta - q^{-1}b^\star\otimes b = \Delta(a^\star).
\end{equation*}
Other claims are verified similarly. Finally, we see that $S(S(\alpha^\star)^\star) = -\alpha$. Since $S(\delta)=-q^{1/2}\gamma$ and $\gamma^\star=q^{-1/2}\beta$, we have
\begin{equation*}
S(S(\alpha^\star)^\star) = S(S(q^{1/2}\delta)^\star) = -q S(\gamma^\star) = -q^{1/2} S(\beta) = - q^{1/2} (q^{-1/2} \alpha) = -\alpha,
\end{equation*}
as expected.

\begin{proposition} 
The algebra ${\cal O}(SP_q^{1|2})$ is left $\star$-quantum superspace for the super Hopf $\star$-algebra ${\cal O}(UOSP_q(1|2))$.
\end{proposition}

\noindent{\it Proof}
We will only see the $\star$ operation for $\xi$. Others can be provided in a similar manner. Since $\delta_L(\xi)=a\otimes \xi + \alpha\otimes x+b\otimes\eta$, one has
\begin{eqnarray*}
[\delta_L(\xi)]^\star &=& a^\star\otimes \xi^\star + \alpha^\star\otimes x^\star + b^\star\otimes \eta^\star \\
&=& q^{1/2} a^\star\otimes \eta + \alpha^\star\otimes x - q^{-1/2} b^\star\otimes \xi \\
&=& q^{1/2} \delta_L(\eta) = \delta_L(\xi^\star)
\end{eqnarray*}
with the first equality in (\ref{11}).

\noindent{\bf Remark.} The algebra ${\cal O}(SP_q^{1|2})$ is not right $\star$-quantum superspace for the super-Hopf $\star$-algebra ${\cal O}(UOSP_q(1|2))$. However, if we define a linear homomorphism $\delta'_R$ on the generators of ${\cal O}(SP_q^{1|2})$ by
$\delta'_R(\xi)=\xi\otimes a - x\otimes \alpha + \eta\otimes b$, {\it etc.} (or in compact form $\delta'_R(X)=X\otimes T^{st}$) we have
\begin{eqnarray*}
[\delta'_R(\xi)]^\star &=& \xi^\star\otimes a^\star - x^\star\otimes \alpha^\star + \eta^\star\otimes b^\star \\
&=& q^{1/2} \eta\otimes a^\star - x\otimes \alpha^\star - q^{-1/2} \xi\otimes b^\star \\
&=& q^{1/2} \delta'_R(\eta) = \delta'_R(\xi^\star).
\end{eqnarray*}

\section{Covariant Differential Calculi on the Superspace $SP_q^{1|2}$}\label{sec5}

In this section, we set up two covariant differential calculi on the symplectic superspace $SP_q^{1|2}$. They contain functions on $SP_q^{1|2}$, their differentials and differential forms.

\subsection{$Z_2$-Graded Differential Algebra}\label{sec5.1}

Let us begin with the definition of the $Z_2$-graded differential calculus. Let $\cal A$ be an arbitrary algebra with unity and $\Gamma$ be a bimodule over
$\cal A$.

\begin{definition}
A first order Z$_2$-graded differential calculus over $\cal A$ is a pair $(\Gamma,{\sf d})$ where ${\sf d}: {\cal A}\longrightarrow\Gamma$ is a linear mapping such that
\begin{equation*}
{\sf d}(fg) = ({\sf d}f)\, g + (-1)^{p(f)} f\, ({\sf d}g) \quad \mbox{for any} \quad f,g\in {\cal A}
\end{equation*}
and where $\Gamma$ is the linear span of elements of the form $a\cdot{\sf d}b\cdot c$ with $a,b,c\in \cal A$.
\end{definition}

\noindent For the product by elements $a\in \cal A$ we will write $aw$ and $wa$, $w\in \Gamma$. The condition in Definition 7 is called the {\it $Z_2$-graded Leibniz rule}.

A $Z_2$-graded differential algebra over $\cal A$ is a $Z_2$-graded algebra $\Gamma=\bigoplus_{n\ge0} \Gamma^n$, $\Gamma^0:= \cal A$, with the linear map
${\sf d}$ of grade one such that ${\sf d}^2 := {\sf d}\circ {\sf d} = 0$ and $Z_2$-graded Leibniz rule holds.

\subsection{Covariance of the Calculus}\label{sec5.2}

Let $\cal A$ be a Hopf algebra, ${\cal X}$ be a quantum space for $\cal A$ with action $\delta_L$ and $\tau: \cal A \longrightarrow \cal A$ is the linear map of degree zero defined by $\tau(\vartheta)=(-1)^{p(\vartheta)}\vartheta$ for all $\vartheta\in \cal A$.

\begin{definition} 
A $Z_2$-graded differential calculus $(\Gamma,{\sf d})$ over ${\cal X}$ is said to be {\sf left (right)-covariant} with respect to $\cal A$ if there exists an algebra homomorphism  $\Delta_L: \Gamma\longrightarrow {\cal A}\otimes \Gamma$ \,$(\Delta_R: \Gamma\longrightarrow \Gamma\otimes \cal A)$ which is a left coaction of $\cal A$ on $\Gamma$  such that:

\noindent(i) $\Delta_L(u)=\delta_L(u)$ \quad $(\Delta_R(u)=\delta_R(u))$ for $u\in {\cal X}$,

\noindent(ii) $\Delta_L\circ{\sf d}=(\tau\otimes{\sf d})\circ\Delta_L$ \quad $(\Delta_R\circ{\sf d}=({\sf d}\otimes\mbox{id})\circ\Delta_R)$.
\end{definition}

\subsection{Commutation Relations} \label{sec5.3}

There exist two covariant $Z_2$-graded first order differential calculi $\Gamma_\pm$ over (1+2)-quantum symplectic superspace $SP_q^{1|2}$ with respect to
the super-Hopf algebra ${\cal A}:={\cal O}(SP_q(1|2))$. One of both is given in the following proposition. The other will be expressed as a result.

\begin{proposition} 
There exists a left covariant $Z_2$-graded first order differential calculus $\Gamma_+$ over $SP_q^{1|2}$ with respect to the super-Hopf algebra ${\cal A}$ such that the set $\{{\sf d} x,{\sf d}\xi,{\sf d}\eta\}$ is a free right $SP_q^{1|2}$-module basis of $\Gamma_+$ and
\begin{eqnarray} \label{13}
x \, {\sf d}x &=& q \, {\sf d}x \, x + q^{1/2}(q^2-1) \, {\sf d}\eta \, \xi, \quad x \, {\sf d}\xi = q \,{\sf d}\xi \,x + (q^2-1) \, {\sf d}x \,\xi,\nonumber \\
x \, {\sf d}\eta & = & q \, {\sf d}\eta \, x, \quad \xi \, {\sf d}x = - q \, {\sf d}x \, \xi, \quad \xi \, {\sf d}\xi = {\sf d}\xi \, \xi, \nonumber \\
\xi \, {\sf d}\eta &=& q^2 \, {\sf d}\eta \, \xi, \quad \eta \, {\sf d}x =- q \, {\sf d}x \, \eta + (1-q^2) \, {\sf d}\eta \, x,\\
\eta \, {\sf d}\xi & = & q^2 \, {\sf d}\xi \, \eta + (1-q^2)[(1+q) \, {\sf d}\eta \, \xi - q^{1/2} \, {\sf d}x \, x], \quad
\eta \, {\sf d}\eta = {\sf d}\eta \, \eta.\nonumber
\end{eqnarray}
\end{proposition}

\noindent{\it Sketch of proof}. To obtain nine cross-commutation relations between the elements of the set $\{\xi,x,\eta\}$ and the elements of the set
$\{{\sf d}\xi,{\sf d}x,{\sf d}\eta\}$, we write, for example, the express $x \,{\sf d}x$ in terms of linearly ${\sf d}x \,x$, ${\sf d}x \,\xi$, ${\sf d}x \,\eta$, {\it etc.} Totally, we have 81 indeterminate coefficients which are successively by the rest of the conditions. Applying the differential {\sf d} to the cross-commutation relations will give eliminate about half of the coefficients. The fact that compatibility with the coaction of ${\cal O}(SP_q(1|2))$ leaves one free parameter.

\begin{corollary} 
The corresponding relations of the second $Z_2$-graded first order differential calculus $\Gamma_-$ over $SP_q^{1|2}$ are obtained if we replace $q$ by $q^{-1}$, $x$ by $\sqrt{-1} \,x$ and $\xi$ by $\eta$ and differentials in the relations (\ref{13}).
\end{corollary}

The relations (\ref{13}) may be written in a compact form
\begin{equation*}
(-1)^{p(X)} X\otimes {\sf d}X = \hat{B_q} {\sf d}X \otimes X.
\end{equation*}
Here the non-zero elements of the matrix $\hat{B_q}$ are
\begin{eqnarray*}
B_{11}^{11} = B^{33}_{33} = -1, \quad B^{12}_{21} = B^{21}_{12} = B^{22}_{22} = B^{23}_{32} = B^{32}_{23} = q, \\
B^{13}_{31} = B^{31}_{13} = -q^2, \quad B^{21}_{21} = B^{32}_{32} = q^2-1,\\
\quad B^{22}_{31} = -B^{31}_{22} = q^{1/2}(q^2-1), \quad B^{31}_{31} = (1+q)(q^2-1).
\end{eqnarray*}
The eigenvalues of the matrix $\hat{B_q}$ are $-1$, $q^2$ and $q^3$ and it admits the spectral decomposition (as a sum of projectors)
\begin{equation*}
\hat{B_q} = -{\sf P}_- + q^2 {\sf P}_+ + q^3 {\sf P}_0
\end{equation*}
where
\begin{eqnarray*}
{\sf P}_- &=& \frac{\hat{B_q}^2 - (q^2+q^3) \hat{B_q} + q^5I}{(q^2+1)(q^3+1)} \\
{\sf P}_+ &=& \frac{\hat{B_q}^2 + (q^3-q) \hat{B_q} + q^3I}{q^2(q-1)(q^2+1)} \\
{\sf P}_0 &=& \frac{\hat{B_q}^2 + (1-q^2) \hat{B_q} - q^2I}{q^2(q-1)(q^3+1)}
\end{eqnarray*}
provided that $q(q-1)(q^2+1)(q^3+1)\ne0$. The projectors obey ${\sf P}_i{\sf P}_j = \delta_{ij}{\sf P}_j$ (no summation) and sum of them equals to the unit matrix. We conclude that the algebra ${\cal O}(SP_q^{1|2})$ has the defining relations ${\sf P}_- X\otimes X = 0$ and that ${\cal O}(SP_q^{1|2})$ is the quotient of $k\{\xi,x,\eta\}$ by the two-sided ideal generated by im${\sf P}_- =$ ker$(\hat{B_q}+I)$.

The proof of the following proposition is obtained by applying the exterior differential {\sf d} to the relations in (\ref{13}).

\begin{proposition} 
The commutation relations between the differentials have the form
\begin{eqnarray} \label{14}
&{\sf d}x\wedge {\sf d}x = q^{1/2}(q-1) \, {\sf d}\eta\wedge{\sf d}\xi, & {\sf d}x\wedge {\sf d}\xi = q^{-1} \, {\sf d}\xi\wedge {\sf d}x, \nonumber \\
&{\sf d}x\wedge {\sf d}\eta = q \, {\sf d}\eta\wedge {\sf d}x, & {\sf d}\xi\wedge {\sf d}\eta = q^2 \, {\sf d}\eta\wedge {\sf d}\xi.
\end{eqnarray}
\end{proposition}

The relations (\ref{14}) can be rewritten as $(-1)^{p({\sf d}X)} {\sf P}_+ {\sf d}X\otimes {\sf d}X = 0$.

\subsection{The Relations With Partial Derivatives}\label{sec5.4}

The calculus will be completed by giving by the following three propositions. We first introduce commutation relations between the coordinates of the quantum superspace and their partial derivatives.

\begin{proposition} 
The relations between the generators of ${\cal O}(SP_q^{1|2})$ and partial derivatives are as follows
\begin{eqnarray} \label{15}
\partial_x x & = & 1 + q \, x \partial_x + (q^2-1) \xi \partial_{\xi}, \quad \partial_x \xi = q \, \xi \partial_x, \nonumber \\
\partial_x {\eta} & = & q \, {\eta} \partial_x + q^{1/2}(1-q^2) \, x\partial_{\xi}, \nonumber \\
\partial_{\xi} x & = & q x \partial_{\xi}, \quad \partial_{\xi} \xi = 1 - \xi \partial_{\xi}, \quad \partial_{\xi} {\eta} = -q^2 {\eta} \partial_{\xi},  \\
\partial_{\eta} x & = & q x \partial_{\eta} + q^{1/2}(q-1) \xi \partial_x, \quad \partial_{\eta} \xi = - q^2 \, \xi \partial_{\eta}, \nonumber \\ \partial_{\eta} {\eta} & = & 1 - {\eta} \partial_{\eta} + (q^2-1) [x \partial_x + (q+1) \, \xi \partial_{\xi}].\nonumber
\end{eqnarray}
\end{proposition}

\noindent{\it Proof}
If $\Gamma$ is the left-covariant differential calculus in Proposition 4, for any element $f$ in ${\cal O}(SP_q^{1|2})$ there are uniquely determined elements $\partial_i(f)$ in ${\cal O}(SP_q^{1|2})$ such that ${\sf d}f = \sum_i {\sf d}x_i \partial_i(f)$. If we replace $f$ with $x_if$ in the left and right hand side of this equality and use the relations (\ref{13}), in results the coefficients of ${\sf d}x_i$ on both sides must coincide. This gives the relations in (\ref{15}).

The proof of the following proposition can be made by using the fact that ${\sf d}^2 = 0$.

\begin{proposition} 
The commutation relations among  the partial derivatives are as follows
\begin{eqnarray} \label{16}
\partial_x \partial_{\xi} &=& q^{-1} \partial_{\xi} \partial_x, \quad \partial_x \partial_{\eta} = q \,\partial_{\eta} \partial_x,
\quad \partial_{\xi}^2=0=\partial_{\eta}^2, \nonumber\\
\partial_{\xi} \partial_{\eta} &=& - q^2 \partial_{\eta} \partial_{\xi} + q^{1/2}(q-1) \, \partial_x^2.
\end{eqnarray}
\end{proposition}

\noindent{\bf Note.} Although the quantum supersphere $r$ is the central element of the algebra ${\cal O}(SP_q^{1|2})$, it satisfies the following relations
with differentials and partial derivatives:
\begin{eqnarray*}
{\sf d}x_i \cdot r &=& q^{-2} r\cdot {\sf d}x_i, \quad x_i\in\{\xi,x,\eta\}\\
\partial_x r &=& q^2 r \partial_x + (1+q)(1-q+q^2) \,x, \\
\partial_\xi r &=& q^2 r \partial_\xi + q^{-1/2}(1+q^3) \,\eta, \\
\partial_\eta r &=& q^2 r \partial_\eta + q^{-1/2}(1+q^3) \,\xi.
\end{eqnarray*}

\section{Symplectic Quantum Weyl Superalgebra}\label{sec6}

In this section, we will define the star operation, denoted by $*$, for partial derivatives and give some related topics. Throughout this section, we will assume that the deformation parameter $q$ is a non-zero complex number such that $\bar{q}=q^{-1}$.

Let ${\cal X}:={\cal O}(SP_q^{1|2})$. In the limit $q\to1$ the relations (\ref{6}) reduce to the defining relations of the superalgebra  ${\cal O}({\mathbb R}^{1|2})$.

Let ${\cal D}$ be the unital algebra with the generators $\partial_\xi$, $\partial_x$, $\partial_\eta$ and defining relations (\ref{16}). It can be shown that the set of monomials
\begin{equation*}
\{\partial_x^k \partial_\xi^l \partial_\eta^m: \, k\in{\mathbb N}_0, l,m=0,1\}
\end{equation*}
forms a vector space basis of ${\cal D}$. In the case $q=1$ the relations (\ref{15}) and (\ref{16}) together with (\ref{6}) reduce to the defining relations of the Weyl superalgebra denoted by ${\cal W}(1|2)$. This motivates following definition.

\begin{definition} 
The quantum symplectic Weyl superalgebra ${\cal W}_q(1|2)$ is the unital algebra generated by $x$, $\xi$, $\eta$ and $\partial_x$, $\partial_\xi$, $\partial_\eta$ and the relations (\ref{6}), (\ref{15}) and (\ref{16}).
\end{definition}

Since the linear map of ${\cal X}\otimes {\cal D}$ to ${\cal W}_q(1|2)$ defined by $u\otimes v \mapsto u\cdot v$ is a vector space isomorphism, one can consider ${\cal X}$ and ${\cal D}$ as subalgebras of ${\cal W}_q(1|2)$ and the set of monomials $\{x^k\xi^l\eta^m\partial_x^{k'} \partial_\xi^{l'} \partial_\eta^{m'}: \, k\in{\mathbb N}_0, l,m=0,1\}$ is a vector space basis of ${\cal W}_q(1|2)$.

The definition of the star operation for the Grassmann generators can be given as follows.

\begin{definition}
If $\theta_1$ and $\theta_2$ are Grassmann generators, the star operation, denoted by $*$, is defined by
\begin{equation*}
(\theta_1\theta_2)^* = \theta_2^*\theta_1^*, \quad (\theta_i^*)^* = \theta_i.
\end{equation*}
\end{definition}

A similar way to the proof of Proposition 1 can be used for the following proposition.

\begin{proposition}
The algebra ${\cal O}(SP_q^{1|2})$ equipped with the involutions determined by
\begin{equation} \label{17}
\xi^* = \xi, \quad x^* = {\bf i} \,x, \quad \eta^* = \eta
\end{equation}
becomes a $*$-algebra.
\end{proposition}

By Proposition 8, the first order differential calculus $\Gamma_+$ given in section 5 is a $*$-calculus for the $*$-algebra ${\cal O}(SP_q^{1|2})$. The involution of $\Gamma_+$ induces the involution for the partial differential operators:

\begin{proposition} 
The algebra ${\cal D}$ supplied with the involutions determined by
\begin{eqnarray} \label{18}
\partial_\xi^* = \partial_\xi, \quad \partial_x^* = q^{-1} \,{\bf i} \, \partial_x, \quad \partial_\eta^* = q^{-2} \partial_\eta,
\end{eqnarray}
becomes a $*$-algebra.
\end{proposition}
The proof of above proposition can be made in a similar way to the proof of Proposition 1.

The involutions given in (\ref{17}) and (\ref{18}) allow us to define the hermitean operators
\begin{equation*}
\hat{x} = (1+{\bf i}) x, \quad \hat{\xi} = \xi, \quad \hat{\eta} = \eta,
\end{equation*}
and
\begin{equation*}
\hat{p}_\xi = \partial_\xi, \quad \hat{p}_x = (q^{-1}+{\bf i}) \partial_x, \quad \hat{p}_\eta = (1+q^{-2})\partial_\eta.
\end{equation*}
Using (\ref{6}), (\ref{15}) and (\ref{16}), new relations will be provided by these operators can be easily obtained. These operators with the resulting relations form a symplectic superalgebra so called the symplectic Heisenberg superalgebra. We will not go into further details. In section 8, we will give a detailed discussion for $h$-deformation.

\section{An Aspect to the $h$-Deformation}\label{sec7}

In this section, we introduce an $h$-deformation of the superspace $SP^{1|2}$ from the $q$-deformation via a contraction following the method of \cite{agha}. Here we denote $q$-deformed objects by primed quantities. Unprimed quantities represent transformed coordinates.

We consider the $q$-deformed algebra of functions on the quantum superspace $SP_q^{1|2}$ generated by $\xi'$, $x'$ and $\eta'$ with the relations (\ref{6}) and we introduce new coordinates $\xi$, $x$ and $\eta$ with the change of basis in the coordinates of the q-superspace using the following matrix $g$:
\begin{eqnarray} \label{19}
X' = \begin{pmatrix} \xi' \\ x' \\ \eta' \end{pmatrix} = \begin{pmatrix} 1 & 0 & \frac{h}{q-1} \\ 0 & 1 & 0 \\ 0 & 0 & 1 \end{pmatrix} \begin{pmatrix} \xi \\ x \\ \eta \end{pmatrix} = g \,X
\end{eqnarray}
where $h$ is a {\it new} deformation parameter that will be replaced by $q$ in the limit $q\to1$.

After the relations (\ref{6}) are used, by taking the limit $q\to1$ we obtain the following exchange relations, which define the symplectic $h$-superspace $SP_h^{1|2}$:

\begin{definition} 
Let ${\cal O}(SP_h^{1|2})$ be the algebra with the generators $\xi$, $x$ and $\eta$ satisfying the relations
\begin{eqnarray} \label{20}
x\xi = \xi x + 2h \eta x, \quad x\eta = \eta x, \quad \xi\eta = -\eta\xi, \quad \xi^2 = h \,(x^2+2\xi\eta), \,\,\, \eta^2 = 0
\end{eqnarray}
where the coordinate $x$ is even and the coordinates $\xi$ and $\eta$ are odd. We call ${\cal O}(SP_h^{1|2})$ the algebra of functions on the $Z_2$-graded quantum symplectic space $SP_h^{1|2}$.
\end{definition}

In the case of dual (exterior) $h$-superspace, denoted by $SP_h^{*1|2}$, we use the transformation
\begin{equation*}
\hat{X}'=g \hat{X}
\end{equation*}
with the components $\theta'$, $y'$ and $z'$ of $\hat{X}'$. The definition is given below.

\begin{definition} 
Let  ${\cal O}(SP_h^{*1|2}):=\Lambda(SP_h^{1|2})$ be the algebra with the generators $\theta$, $y$ and $z$ satisfying the relations
\begin{equation*}
\theta^2 = h z^2, \quad \theta y = y \theta - 2h z\theta, \quad \theta z = z \theta, \quad yz = zy + 2 h z^2
\end{equation*}
where the coordinate $\theta$ is odd and the coordinates $y$ and $z$ are even. We call $\Lambda(SP_h^{1|2})$ the quantum dual (exterior) algebra of the $Z_2$-graded quantum space $SP_h^{1|2}$.
\end{definition}

Obviously, in the limit $h\to0$ the algebra ${\cal O}(SP_h^{1|2})$ is the $Z_2$-graded polynomial algebra in three supercommuting indeterminates and the algebra $\Lambda(SP_h^{1|2})$ is the exterior algebra of $SP_h^{1|2}$.

\begin{definition}[\cite{kac}]
A Lie superalgebra $L$ is a superalgebra $L = L_{\bar{0}}\oplus L_{\bar{1}}$ with an operation $[,]$ satisfying the following axioms: for $a\in L_\alpha$, $b\in L_\beta$
\begin{eqnarray*}
\left[a, b\right] &=& -(-1)^{\alpha.\beta} [b, a] \qquad (\mbox{supersymmetry})\\
\left[a, [b, c]\right] &=& [[a, b], c] + (-1)^{\alpha.\beta} [b, [a, c]] \qquad (\mbox{Jacobi identity}).
\end{eqnarray*}
\end{definition}

Following proposition can be easily proven by using the Definition 13.
\begin{proposition} 
The superalgebras ${\cal O}(SP_h^{1|2})$ and $\Lambda(SP_h^{1|2})$ are both Lie superalgebra.
\end{proposition}

The $h$-deformed matrix $C$ and the matrix $\hat{B}$ are given by
\begin{equation} \label{21}
C_h = \lim_{q\to1} \left(g^{st}C_qg\right) = \begin{pmatrix} 0 & 0 & -1 \\ 0 & 1 & 0 \\ 1 & 0 & h \end{pmatrix}, \quad \hat{B_h} = \lim_{q\to1} \left[(g\otimes g)^{-1}\hat{B_q}(g\otimes g)\right],
\end{equation}
respectively. Here the matrix $C_h$ corresponds to the matrix $(C^{-1})^t$ in \cite{Jus-Sob} with $h=p$, the matrix $B_h=P\hat{B_h}$ corresponds to the matrix $B$ in \cite{Aiz-Chak3} with replacing $h$ to $-h/2$. Using the $h$-deformed symmetrizer $P_{h-}$,
\begin{equation*}
P_{h-} = \lim_{q\to1} \left[(g\otimes g)^{-1} P_- (g\otimes g)\right],
\end{equation*}
we can rewrite the relations in (\ref{20}) in a compact form ${\sf P}_{h-} X\otimes X = 0$.

As an interesting case, ($h$-deformed) supersphere on the symplectic $h$-super\-space is given by
\begin{equation*}
\rho = X^{st}C_hX = x^2 + 2\xi\eta = h^{-1}\xi^2
\end{equation*}
where $h\ne0$.

The corresponding $h$-deformation of the supergroup $SP(1|2)$ as a matrix quantum supergroup $SP_h(1|2)$ generated by the matrix elements of the matrix $T$ in (\ref{1}) can be obtained from the requirement that $SP_h^{1|2}$ and $SP_h^{*1|2}$ have to be covariant under the left coaction
\begin{equation} \label{22}
\delta_L(X) = T \dot{\otimes} X, \quad \delta_L(\hat{X}) = T \dot{\otimes} \hat{X}
\end{equation}
and assuming that the matrix elements of the matrix $T$ supercommute with the elements of $SP_h^{1|2}$ and $SP_h^{*1|2}$. Now the orthosymplectic condition
\begin{equation*}
T^{st}C_hT = D_hC_h, \quad TC_h^{-1}T^{st} = D_hC_h^{-1}
\end{equation*}
and (\ref{22}) determine the relations among the entries of $T$ similar to the relations in \cite{Aiz-Chak3}.

Relations describing $h$-deformation of the quantum supergroup $SP_h(1|2)$ can be alternatively obtained by performing the similarity transformation introduced in \cite{agha}
\begin{equation*}
T' = g \, T \, g^{-1}.
\end{equation*}
Here $T'$ is a quantum super matrix in $SP_q(1|2)$ and then the matrix elements of $T'$ satisfy the relations (\ref{4}).

\section{Symplectic $h$-Deformed Heisenberg Superalgebra}\label{sec8}

Using some relations in previous sections, we can obtain a $h$-deformed symplectic Heisenberg superalgebra. Let us begin the following definition.

To realize the superalgebras, one has to use not only the even but also the odd operators \cite{kac}. So, it is of importance to work the Heisenberg superalgebra denoted by ${\cal H}$. This algebra generated by even coordinate $x$, Grassmann coordinates $\xi$, $\eta$, and even operator $p_x$, Grassmann operators $p_\xi$, $p_\eta$. These elements of this algebra are supposed to be selfadjoint.

The coordinates and operators satisfy the super-commutation relations with the additional relations
\begin{eqnarray*}
\left[p_x, x\right] &=& 2{\bf i}{\bf 1}, \quad [p_x, \xi] = \left[p_x, \eta\right] = 0, \nonumber\\
\left[p_\xi, x\right] &=& 0, \quad [p_\xi, \xi]_+ = {\bf 1}, \quad [p_\xi, \eta]_+ = 0, \\
\left[p_\eta, x\right] &=& 0, \quad [p_\eta, \xi]_+ = 0, \quad \left[p_\eta, \eta\right]_+ = {\bf 1}.\nonumber
\end{eqnarray*}

We wish to modify the superalgebra ${\cal H}$ according to symplectic quantum supergroup ideas. For this, for example, it is natural to assume that
$\hat{x}$ is a coordinate of a quantum supervector and to relate $\hat{p}_x$ to a derivative in such a supervector. To identify the momenta ${\bf i} \,p_x$, $p_\xi$ and $p_\eta$ with $\partial_x$, $\partial_\xi$ and $\partial_\eta$, we must take care of hermiticity of coordinates and momenta. Let us begin with the following proposition.

\begin{proposition} 
The coordinate functions with the partial derivatives satisfy the following $h$-deformed relations
\begin{equation} \label{23}
\partial_i X^j = \delta_{ij} + \sum_{k,l} \hat{B}^{jk}_{il} X^l \partial_k, \quad \sum_{k,l} ({\sf P}_{h-})^{lk}_{ij} \,\partial_k \partial_l = 0.
\end{equation}
\end{proposition}

\begin{definition} 
We assume that $h$ is a complex number such that $\bar{h}=-h$. Then the hermitean conjugation of the coordinates $\xi$, $x$ and $\eta$, and the partial derivatives $\partial_x$, $\partial_\xi$ and $\partial_\eta$ can be defined by
\begin{equation} \label{24}
\xi^* = \xi - h\eta, \quad x^* = {\bf i} \, x, \quad \eta^* = \eta,
\end{equation}
and
\begin{equation} \label{25}
\partial_\xi^* = \partial_\xi, \quad \partial_x^* = {\bf i} \, \partial_x, \quad \partial_\eta^* = \partial_\eta - h \partial_\xi,
\end{equation}
respectively.
\end{definition}
\noindent The conjugation, thus, defined is an involution because its square is the identity.

The relations (\ref{6}), (\ref{15}) and (\ref{16}) are now invariant under the definitions in (\ref{24}) and (\ref{25}). The involution above allows us to define the hermitean operators
\begin{equation*}
\hat{x} = (1+{\bf i}) x, \quad \hat{\xi} = \xi, \quad \hat{\eta} = \eta + h x,
\end{equation*}
and
\begin{equation*}
\hat{p}_\xi = \partial_\xi, \quad \hat{p}_x = (q^{-1}+{\bf i}) \partial_x, \quad \hat{p}_\eta = (1+q^{-2})\partial_\eta.
\end{equation*}
Then, the following proposition can be proven using the relations in (\ref{23}) and Definition 13.

\begin{proposition} 
(i) The $h$-deformed symplectic Heisenberg superalgebra ${\cal H}_h$ is generated by the elements of the set
$\{\hat{p}_x,\hat{p}_\xi,\hat{p}_\eta,\hat{x},\hat{\xi},\hat{\eta},{\bf 1}\}$ with the following commutation relations
\begin{eqnarray*}
\left[\hat{x}, \hat{\xi}\right] &=& 2h\hat{\eta}\hat{x}, \quad \left[\hat{x}, \hat{\eta}\right] = 0, \quad [\hat{\eta}, \hat{\xi}]_+ = 0, \\
\hat{\xi}^2 &=& -h(\mbox{$\frac{{\bf i}}{2}$} \, \hat{x}^2 - 2\hat{\xi}\hat{\eta}), \quad \hat{\eta}^2 = 0, \\
\left[\hat{p}_x, \hat{p}_\xi\right] &=& 0, \quad \left[\hat{p}_x, \hat{p}_\eta\right] = -2h\hat{p}_\xi\hat{p}_x, \quad \left[\hat{p}_\xi, \hat{p}_\eta\right]_+ = 0, \\
\hat{p}_\xi^2 &=& 0, \quad \hat{p}_\eta^2 = -h(\mbox{$\frac{{\bf i}}{2}$} \, \hat{p}_x^2 + \hat{p}_\eta\hat{p}_\xi), \nonumber\\
\left[\hat{p}_x, \hat{x}\right] &=& 2{\bf i} \,({\bf 1} + 2h \hat{\eta} \hat{p}_\xi), \quad [\hat{p}_x, \hat{\xi}] = 2h \hat{x} \hat{p}_\xi, \quad \left[\hat{p}_x, \hat{\eta}\right] = 0, \quad \left[\hat{p}_\xi, \hat{x}\right] = 0,\\
\left[\hat{p}_\xi, \hat{\xi}\right]_+ &=& {\bf 1}+2h\hat{\eta} \hat{p}_\xi, \quad [\hat{p}_\xi, \hat{\eta}]_+=0, \quad \left[\hat{p}_\eta, \hat{x}\right]= 0,\\
\left[\hat{p}_\eta, \hat{\xi}\right]_+ &=& \mbox{$\frac{h}{2}$} \left(2{\bf i} \hat{x}\hat{p}_x - (\hat{\xi + \hat{\eta}}) \hat{p}_\xi + \hat{\eta} (\hat{p}_\eta - \hat{p}_\xi)\right), \quad \left[\hat{p}_\eta, \hat{\eta}\right] = {\bf 1}+2h\hat{\eta} \hat{p}_\xi.
\end{eqnarray*}

(ii) The  symplectic Heisenberg superalgebra ${\cal H}_h$ is a Lie superalgebra.
\end{proposition}


\baselineskip=10pt
\begin{thebibliography}{}
\bibitem{agha} Aghamohammadi, A, Khorrami, M., Shariati, A.: $h$-deformation as a contraction of $q$-deformation, J. Phys. A: Math. Gen. {\bf 28}, L225-L231 (1995).
\bibitem{Aiz-Chak1} Aizawa, N and Chakrabarti, R.: Quantum Spheres for $OSP_q(1|2)$, J. Math. Phys. {\bf 46}, 103510- (2005).
\bibitem{Aiz-Chak2} Aizawa, N and Chakrabarti, R.: Noncommutative Superspaces Covariant Under $OSP_q(1|2)$ Algebra, Lie Theory and Its Applications in Physics VI ed. V.K. Dobrev et al, Heron Press, Sofia, 2006.
\bibitem{Aiz-Chak3} Aizawa, N and Chakrabarti, R.: Noncommutative geometry of super-jordanian $OSP_h(2|1)$ covariant quantum space, J. Math. Phys. {\bf 45}, 1623-1638 (2004).
\bibitem{Ber-Tos} Berezin, F. A., Tolstoy, V. N.: The group with Grassmann structure $UOSP(1|2)$, Commun. Math. Phys. {\bf 78}, 409-428 (1981).
\bibitem{Celik1} Celik, S.: Differential geometry of the $q$-superplane, J. Phys. A: Math. Gen. {\bf 31}, 9695-9701 (1998).
\bibitem{Connes} Connes, A.: Non-commutative differential geometry, Publ. IHES 62, 257-360 (1985).
\bibitem{FRT} Faddeev, L. D., Reshetikhin N. Yu., and Takhtajan, L. A.: Quantization of Lie groups and Lie algebras, Leningrad Math. J. {\bf 1}, 93-225 (1990).
\bibitem{Fra-Sci-Sor} Frappat, L., Sciarrino, A., and Sorba, P.: Dictionary on Lie superalgebras, hep-th/9607161.
\bibitem{Jus-Sob} Juszczak C. and Sobczyck, J.: New quantum deformation of $OSP_h(1|2)$, Czech. J. Phys. {\bf 48}, 1375-1383 (1998).
\bibitem{kac} Kac, V.: Lie Superalgebras, Adv. in Math. {\bf 26}, 8-96 (1977).
\bibitem{Kul-Res}  Kulish P. P. and Reshetikhin,  N. Yu.: Universal $R$-matrix of the quantum superalgebra $osp(2|1)$, Lett. Math. Phys. {\bf 18}, 143-149 (1989).
\bibitem{Manin1} Manin, Yu I.: Quantum groups and noncommutative geometry, Montreal Univ. Preprint, 1988.
\bibitem{Manin2} Manin, Yu I.: Multiparametric quantum deformation of the general linear supergroup, Commun. Math. Phys. 123, 163-175 (1989).
\bibitem{Rit-Sch} Rittenberg, V., Scheunert, M.: Elementary construction of graded Lie groups, J. Math. Phys. {\bf 19}, 709-713 (1978).
\bibitem{Saleur-Kaufmann} Saleur, H. and Wehefritz-Kaufmann, B.: Integrable quantum field theories with supergroup symmetries: the $OSP(1|2)$ case, Nucl. Phys. B {\bf 663}, 443-466 (2003).
\bibitem{Soni} Soni, S.: Differential calculus on the quantum superplane, J. Phys. A: Math. Gen. {\bf 24}, 619-624 (1990).
\bibitem{Wess-Zumino} Wess J. and Zumino, B.: Covariant differential calculus on the quantum hyperplane, Nucl. Phys. B {\bf 18}, 302-312 (1990).
\bibitem{Woronowicz1} Woronowicz, S.L.: Compact matrix pseudogroups (quantum groups), Comm. Math. Phys. {\bf 111}, 613-665 (1987).
\bibitem{Woronowicz2} Woronowicz, S.L.: Differential calculus on compact matrix pseudogroups (quantum groups), Comm. Math. Phys. {\bf 122}, 125-170 (1989).

\end{thebibliography}
\end{document}